# Gaussian Processes for Observational Dose-Response Inference

Jake Dailey

**Abstract**

We adapt Gaussian processes (GPs) for estimating the average dose-response function (ADRF) in observational settings, introducing a powerful complement to treatment effect estimation for understanding heterogeneous effects. We incorporate samples from a GP posterior for the propensity score into a GP response model using Girard's approach to integrating over uncertainty in training data [12]. We show Girard's method admits a positive-definite kernel, and provide theoretical justification by identifying it with an inner product of kernel mean embeddings. We demonstrate double robustness of our approach under a misspecified response function or propensity score. We characterize and mitigate regularization-induced confounding [15] in GP response models. We show improvement over other methods for ADRF estimation in terms of coverage of the dose-response function and estimation bias, with less sensitivity to misspecification across experiments.

## 1 Introduction

Across Amazon, we are faced with questions of resource allocation that require an understanding of noisy causal relationships. It is not sufficient to study whether or not to allocate a resource; we need to assess *how much* of an input is required to effect some level of output, and when effects diminish.

The majority of observational causal inference literature focuses on estimating average treatment effect (ATE), $E[Y_i(1) - Y_i(0)]$, which estimates the impact of a binary (i.e., treated (1) or not treated (0)) treatment on an individual response $Y_i \in \mathcal{Y} \subset \mathbb{R}$, or conditional average treatment effect (CATE), $E[Y_i(1) - Y_i(0) | \boldsymbol{X}_i = \boldsymbol{x}]$, which specifies to subpopulation $\boldsymbol{x} \in \mathcal{X}$ (for set $\mathcal{X}$, with bolded font denoting a vector throughout). In settings with continuous treatments and continuous outcomes, it is often not appropriate to average effects across treatment levels; instead, we would like to model

$$E[Y_i(t)|T_i = t] = E[E[Y_i(t)| \boldsymbol{X}_i, T_i = t]| T_i = t], \qquad t \in \mathcal{T} \subset \mathbb{R}$$

typically referred to as the average dose-response function (ADRF). Note that here, $Y_i(t)$ is a random function indexed by dosage $t$. In cases where (C)ATE estimates are desired, modeling the ADRF allows us to better-understand heterogeneity in dosage and effects for decision-making. Zhao et al. [37] provide a thorough review of the ADRF literature.

We echo Zhao et al.'s cautionary message, "*estimating the full [A]DRF with a continuous treatment in an observational study is challenging*" [37]. The greatest risk cited is bias due to model misspecification, which we too observe in our experiments. Many popular methods [3, 5, 7, 16, 17] impose a rigid generating process for the response surface by explicitly estimating the contribution of covariates, treatment, and the propensity for treatment. It is well-known that this requires each term to be correctly specified to achieve unbiased causal effect estimates (unless using doubly robust estimation [8]). We eschew these formulations by modeling the ADRF posterior and continuous propensity function as proceeding from GP priors. Rather than explicitly estimating the effects of covariates, treatment, and propensity for treatment, ADRF estimates are then functions of observed responses at nearby points (in a sense we will make precise). Our methods seem to reduce the risk of misspecification; we report a considerable improvement over other ADRF estimation routines.

The primary contributions of this paper are as follows:

- We adapt Gaussian process (GP) regression to the observational dose-response setting by conditioning on a GP-estimated generalized propensity score (PS).
- We prove Girard's modification to the Gaussian kernel [12], which averages over Gaussian-distributed uncertainty in inputs, is itself a positive-definite (p.d.) kernel (and so yields valid GP regression); we term it the Predictive Radial Basis Function (PRBF) kernel.
- We provide theoretical motivation for the PRBF kernel, showing it is an inner product of kernel mean embeddings (KMEs) [13, 29, 30] and that the corresponding kernel matrix contains all information about the Maximum Mean Discrepancy between distributions [13, 30]; we observe in this context that the kernel expresses GP covariance as "concurrence" in KMEs of a GP prior's hypotheses under its distribution-valued arguments.
- We prove *any* inner product of KMEs with respect to a translation-invariant kernel is itself a p.d. kernel, suggesting alternative PRBF kernels with different KMEs are available.
- We characterize regularization-induced confounding (RIC), originally described in [15] in kernel methods, propose the Additive Predictive Radial Basis Function (A-PRBF) kernel and priors for mitigating RIC, and demonstrate a reduction in RIC in experiments.
- We demonstrate double robustness when PS is estimated with a GP regressor.

We first review literature on Gaussian processes and observational causal inference, describe our methods, then apply our method to the simulated dataset described in [16] and a continuous adaptation of the Infant Health and Development Program (IHDP) dataset [14].

## 2 Background

### 2.1 Gaussian processes

GPs flexibly model distributions of beliefs about the relationship between their input and outcome variables; see [23] for a comprehensive introduction. Let $\mathcal{GP}$ be a GP that gives rise to a Gaussian posterior distribution $f$ defined over random vector $\boldsymbol{\theta} \in \mathbb{R}^d$,

$$f(\boldsymbol{\theta}) \sim \mathcal{GP}(m(\boldsymbol{\theta}), k(\boldsymbol{\theta}, \boldsymbol{\theta}')). \tag{1}$$

In the causal inference setting, $\boldsymbol{\theta}$ can represent covariates, treatment variables, PSs, or all of the above; $f$ is the outcome-generating distribution parameterized by $\boldsymbol{\theta}$. $\mathcal{GP}$ (and by extension, $f$) is characterized by a mean function $m$ and kernel function $k$; these can roughly be understood as encoding priors about how the mean of the process varies in $\boldsymbol{\theta}$ and how points in the feature space relate to each other, respectively. When $m = \mathbf{0}$ (as is common but not necessary), we can sample from an updated GP posterior exactly with

$$\hat{Y} \sim f(\Theta_*) = \mathcal{N}(\bar{Y}_*, cov(Y_*)), \text{ where } k(\boldsymbol{\theta}_i, \boldsymbol{\theta}_j) = K_{i,j},$$
$$\bar{Y}_* \triangleq E[f(\Theta_*)|Y] = K(\Theta_*, \Theta)[K(\Theta, \Theta) + \epsilon^2 I]^{-1} Y, \text{ and} \tag{2}$$
$$cov(Y_*) = K(\Theta_*, \Theta_*) - K(\Theta_*, \Theta)[K(\Theta, \Theta) + \epsilon^2 I]^{-1} K(\Theta, \Theta_*),$$

where $i, j$ index units, $\Theta \in \mathbb{R}^{n \times d}$, $\Theta_*$ is a random matrix of new observations, and $\hat{Y}$ is a vector-valued sample from the Gaussian process posterior evaluated at $\Theta_*$. In the ADRF context, GP response inferences can readily be explained in relation to other observations either directly through $k$ (or by noting that they are equivalent to a linear combination of training points weighted by the covariance function, shown in Appendix A.1). GPs allow us to reason about a posterior distribution of possible responses at each value of $\boldsymbol{\theta}$ and construct credible intervals from posterior samples.

GPs have proven useful in many causal inference settings [1, 11, 22, 24, 28, 32, 33, 36]. GP regressors have been used to model dose-response functions in experimental settings [28, 32]. Many authors have proposed adaptations of GP priors for observational ATE estimation [1, 11, 22, 24, 33, 36]. If we observe realizations of a GP, the GP regressor is the best linear unbiased predictor (see Appendix A.1 for result summary). To the author's knowledge, this paper is the first to propose GPs for observational dose-response inference. GPs are a suitable choice for *a priori* smooth and relatively low-dimension ADRF estimations because they allow for straightforward interpretation and quantification of uncertainty. GP convergence is known to slow when the covariance function is specified to be overly-smooth for the data and as the dimension of the feature space increases [31, 34]. This paper finds GPs offer a promising model for studying observational ADRF estimation problems, but require adaptation to address these problems' complexities.

**2.2 Generalized propensity scores**

In observational causal inference, most methods are interested in controlling for the imbalanced application of treatment resulting from self-selection into treatment or other confounding processes. Rosenbaum & Rubin addressed this problem by conditioning a linear outcome model on an estimate of the propensity of individual $i$ to receive a binary treatment [27]. Other authors [17, 18, 26] generalize this approach to continuous treatments by conditioning response models on estimates of the conditional likelihood of a given dose *or* an estimate of the mean dosage. We adopt the latter approach, conditioning on an estimate $\hat{\pi}(x_i)$ of $\pi_i = E[T_i|X_i]$.

Hahn et al. propose such propensity scores to mitigate regularization-induced confounding (RIC) in modeling response surface posteriors [15]. In their characterization, RIC occurs when variation in response caused by *propensity for dosage* is misattributed to the dosage variable, which commonly occurs when regularizing terms induce inappropriate or ambiguous priors over contributions of $\pi$ to the response. Hahn et al. [15] characterize the PS as a covariate-dependent prior which reallocates this variation from $T$ to $\mu$; for example, they model a response to binary treatment as

$$y_i = \mu(x_i, \hat{\pi}(x_i)) + \tau(x_i)t_i + \varepsilon_i, \ \varepsilon_i \sim \mathcal{N}(0, \sigma^2) \quad (3)$$

where $x_i$ and $t_i$ are observed covariates and treatment, $\hat{\pi}(x_i)$ is an estimated propensity score, and the functions $\mu$ and $\tau$ can be non-linear or linear (in their case, they use a forest of Bayesian additive regression trees [5, 16] to model these). See [15], sections 4 and 8 for a complete discussion.

Many authors assume a linear, parametric form for continuous PS models [3, 7, 17, 19, 26]. In the case that $\pi$ is a non-linear function of the covariates, our model will be misspecified and ADRF estimates biased. Many methods explicitly model the respective effects of covariates, treatment, and PSs to estimate the ADRF [16, 36, 37]. In this case, PS and outcome models both need to be correctly specified to achieve unbiased causal effect estimates, which is difficult to ensure in practice. In the ADRF setting, it is unnecessary to isolate each of these effects; the goal is rather to understand how the response surface varies with respect to dosage while controlling for potential confounders. Our approach only requires we correctly specify the contribution of confounders to the covariance (rather than the response). This is useful in the common case where we can reason about appropriate notions of covariance more easily than relationships between covariate, propensity, treatment, and response; experimental results suggest we reduce risk of misspecification.

Along these lines, noting that modeling ADRF as a GP with an inappropriately-smooth covariance or high-dimensional covariates slows convergence [31, 34], conditioning on sufficiently smooth and expressive models for $\pi_i$ can allow us to decrease response model dimensionality by summarizing covariates that only indirectly explain variation in the ADRF via $\pi_i$.

## 3 Method

We first develop a GP model for $\pi_i$, then model the ADRF with a kernel that integrates over $\pi_i$ posterior variance to incorporate this uncertainty into the response posterior. We assume strong ignorability of treatment assignment given covariates, $Y_i(t) \perp\!\!\!\perp T_i | X_i$, and $p(T_i = t | X_i) > 0$ for all $t \in \mathcal{T}$. We caution that these are strong assumptions that will not hold in many practical situations, but are generally considered necessary for unbiased ADRF estimation and motivation for PSs [37].

**3.1 Modeling expected dosage**

We assume that propensity for dosage is unknown and model $\pi(x) = E[T|X = x]$ as realizations of a GP evaluated at $x$ with all observations $x'$,

$$\pi(x) \sim \mathcal{GP}(m_\pi(x), k_\pi(x, x')). \quad (4)$$

$m_\pi$ and $k_\pi$ are selected to reflect analyst beliefs about the dosage-generating process. In some applications, a zero-mean prior may not be appropriate (e.g., the "default" treatment level is some other constant or scales linearly as a function of an observed covariate). In other examples, treatments may vary as non-smooth or periodic functions of training points, informing kernel choice; see [6, 23] for a review of kernel design considerations. Modern machine learning methods (GPs included) are known to carry the risk of biased treatment effect estimates due to overfitting [4], which is mitigated using sample splitting and propensity model selection based on generalization heuristics (e.g., out-of-sample mean squared error) [4]. We employ this approach and retain data efficiency by splitting our dataset in two, first fitting separate PS models on either fold, evaluating those models over their respective holdout sets, then using PSs estimates over the holdout sets in our response model. Samples from both models are pooled to estimate the ADRF. We select $m_\pi$ and $k_\pi$ to minimize out-of-sample bias $\sum_{i=1}^{n} \hat{\pi}(x_i) - T_i$. Optimizing for out-of-sample bias controls sampling variance; minimizing, say, mean squared error can result in high-variance extrapolations in regions of PS posterior with low support by over-emphasizing extreme levels of treatment.

### 3.2 Incorporating $\hat{\pi}$ into a response-model kernel

Girard provides an expression for an RBF kernel averaged over conditionally Gaussian-distributed inputs, describing it as "a Gaussian kernel with its length-scales bounded below" [12, section 3.6]; we provide a derivation in Appendix A.2. The resulting covariance function is given by

$$k_{\hat{\pi}}^{PRBF}(x, x') := \frac{\gamma^2}{\sqrt{2\pi(l^2 + \hat{\sigma}(x)^2 + \hat{\sigma}(x')^2)}} \exp\left(-\frac{(\hat{\mu}(x) - \hat{\mu}(x'))^2}{2(l^2 + \hat{\sigma}(x)^2 + \hat{\sigma}(x')^2)}\right). \qquad (5)$$

where $\hat{\mu}(x)$ is the posterior mean, $\hat{\sigma}(x)^2$ the posterior variance of the propensity score GP, evaluated at $x$, and $l$ is the length-scale parameter. Note we write $k_{\hat{\pi}}^{PRBF}$ because it is characterized by $x \mapsto \hat{\pi}(x) = \mathcal{N}(\hat{\mu}(x), \hat{\sigma}(x))$. It is not obvious that $k_{\hat{\pi}}^{PRBF}$ is a positive-definite kernel (and actually describes a valid GP regressor), nor is it obvious why this way of averaging is appropriate or useful. In Appendix A.2.1, we prove that $k_{\hat{\pi}}^{PRBF}$ is indeed a p.d. kernel and thus term $k_{\pi}^{PRBF}$ the Predictive Radial Basis Function kernel, with "predictive" recalling predictive averaging of kernel values over uncertain inputs. In Appendix A.2.2, we identify $k_{\hat{\pi}}^{PRBF}$ with inner products of kernel mean embeddings, and recognize its kernel matrix $K_{\hat{\pi}}^{PRBF}$ contains all information conveyed by Maximum Mean Discrepancy (MMD) [13, 30]. In the context of zero-mean GP regression, this implies $K_{\pi}^{PRBF}$ summarizes

$$E_{f \sim \mathcal{GP}(0,k)}[(\varphi_x f - \psi_{x'} f)^2] = var[\varphi_x f - \psi_{x'} f],$$

which we refer to as the "concurrence" between Gaussian distributions $\varphi_x$ and $\psi_{x'}$'s embeddings, with expectations taken over hypotheses $f$ with respect to a GP prior with Gaussian kernel $k$. This interpretation of $k_{\hat{\pi}}^{PRBF}$ motivates our approach, and suggests other generalizations under suitable families for $k$, $\varphi_x$, and $\psi_{x'}$ may be available. Some authors note that response models do not require estimates of uncertainty in the PS for valid inference [15, 24]. Others observe point estimates of treatment propensity often have high variance and so magnify the impact of model misspecification [36, 37]. Indeed, inference can be made more robust to misspecification insofar as dosage variation explained by unobserved covariates can be modeled as variance in dosage. $k_{\hat{\pi}}^{PRBF}$'s lower bound on length-scale discounts deviations in dosage according to expected variation in dosage at that point.

### 3.3 Specifying the response model

Hahn et al. demonstrate the importance of separate, pointwise control of the variation of an outcome function $f$ in $T$ to discovering heterogeneous effects [15]. To accomplish this, we define

$$k_{A \cdot PRFB}(\theta, \theta') := k_{\hat{\pi}}^{PRBF}\big((\hat{\pi}(x), x), (\hat{\pi}(x'), x')\big) + k_{RBF}(t, t'), \text{ and} \qquad (6)$$

$$k_{RBF}(t, t') = \omega^2 exp\left(-\frac{(t - t')^2}{2\rho^2}\right),$$

where $\theta = (\pi(x), x, t)$ is the vector formed by concatenating covariate and treatment variables. We abuse notation in writing $k_{\hat{\pi}}^{PRBF}\big((\hat{\pi}(x), x), (\hat{\pi}(x'), x')\big)$ to denote a PRBF kernel evaluated over PS posterior distributions *and* covariates; $\hat{\sigma}(\cdot)^2 = 0$ for (non-PS) covariate dimensions without uncertainty estimates in $k_{\hat{\pi}}^{PRBF}$. This additive kernel corresponds to

$$f\big((\hat{\pi}(x), x)\big) + g(t) = h(\theta) \sim \mathcal{GP}\big(m(\theta), k_{A \cdot PRFB}(\theta, \theta')\big), \qquad (7)$$

with $f \sim \mathcal{GP}(0, k_{\hat{\pi}}^{PRBF})$ and $g \sim \mathcal{GP}(0, k_{RBF})$. We show this specification is doubly robust to misspecification of either the response function ($f$ and $g$) or the propensity score ($\hat{\pi}$) in Appendix A.3. $k_{A \cdot PRFB}$ mitigates a different form of RIC than the exact phenomenon described by Hahn et al. [15]. By default, maximizing the log marginal likelihood of the response with respect to a non-additive model's parameters risks misattributing likelihood contributions due to *covariation* along $\pi$ to covariation along $t$. In other words, vague priors for mean and kernel functions may have the unintended effect of over-attributing the influence of $t$ on the response in cases where influence should have been attributed to $\pi_i$. Additivity allows us to explicitly and independently regularize the contribution of $t$ to covariance, relative to $\pi$ and $x$. An alternative mitigation is achieved by allowing $l$ to be a vector and placing a more conservative prior over $t$'s dimension relative to $\hat{\pi}$; to some extent, $k_{\hat{\pi}}^{PRBF}$ implicitly does this by lower-bounding $l$ along $\hat{\pi}$. Modeling a zero-mean GP equipped with $k_{A \cdot PRFB}$ requires only that we correctly specify $k_{\hat{\pi}}^{PRBF}$'s contribution to the covariance rather than its contribution to the response; it is robust to feature selection for $x$ and effects the same posterior insofar as $k_{\hat{\pi}}^{PRBF}$'s evaluation captures the additive contribution covariance of confounders.

The response GP is specified per (7), with $m$ selected based on practitioner beliefs about the mean of the response process. Given response observations $y$, we adapt Hahn et al. [15] and Gelman [10]'s approaches, applying half-Normal priors to both sub-model's scale parameters such that the half-Normal median is *twice* the marginal standard deviation of $y$ for $\gamma$ and *half* the marginal standard deviation of $y$ for $\omega$. These place non-trivial density around 0 but reflect the belief that posteriors are primarily determined by variation in $\hat{\pi}$ rather than $t$. The more we expect the response to be influenced by covariance along $\hat{\pi}$ rather than $t$ itself, the smaller $\omega$ should be relative to $\gamma$. We apply half-Normal priors over $l$ and $\rho$ based on scale of each dimension. We update the GP prior by maximizing the log marginal likelihood (LML) of $y$ with respect to $\gamma$, $\omega$, $\rho$, and $l$, then sample DRF estimates from the GP posterior conditioned on observed covariates, treatments, and estimated PSs.

## 4 Empirical evaluations

We test our model in fully-simulated and semi-simulated scenarios. In all studies, we estimate the 90% interval coverage ($Cov_{90}$), mean 90% interval length ($I_{90}$), bias, and root mean squared error (RMSE) of the ADRF. $Cov_{90}$ is the proportion of observations whose interval includes the true ADRF at their treatment level $t$. $I_{90}$ is defined as the average difference in interval bounds. Since we are estimating continuous functions rather than a point estimate, "good" ADRF models should have high coverage and *acceptable* interval lengths, biases, RMSEs; the latter two metrics in particular can mislead by averaging across localized performance.

For GP models, bias and RMSE are averaged over all units' posterior-sampled errors and are relative to the true ADRF at their treatment level $t$. For GP models, 90% intervals are credible intervals over the response model posterior. For non-GP models (except BART), we bootstrap sample ADRF estimates 50 times; we construct 90% intervals using quantiles over sampled ADRF estimates, and average bias and RMSE over bootstrapped, unit-level errors. See Appendix D.1 for formal definitions of evaluation metrics.

### 4.1 Propensity score GP specifications

We implement an exact GP in GPyTorch [9] with zero mean and a 1/2-Matern kernel equipped with a scaling parameter. We iteratively maximize the LML of $t$ with respect to the parameters to convergence (after 1000 epochs using the Adam optimizer and a learning rate of 0.0015). All GPs are regularized by including a zero-mean Gaussian noise term, as described in [23].

### 4.2 Response GP specifications

All GP priors assume zero mean. For all experiments, we maximize the LML of $y$ with respect to kernel parameters to convergence (after 7000 epochs using Adam optimizer with a learning rate of 0.0025). All GPs are regularized by including a zero-mean Gaussian noise term, as in [23]. We do an ablative comparison of our proposed kernel ($k_{A \cdot PRFB}$, as in Section 3):

Additive Radial Basis Function (A-RBF): This kernel is identical to $k_{A \cdot PRFB}$, but we set all estimated $\pi$ posterior variances to 0, effectively reverting the first term of $k_{A \cdot PRFB}$ to a standard RBF over covariates and treatment variables. This variation allows us to isolate the effect of modeling uncertainty in $\pi$ as part of an additive kernel via $k_{\hat{\pi}}^{PRBF}$.

Predictive Radial Basis Function (PRBF): This is exactly $k_{\hat{\pi}}^{PRBF}$ as described in 3.2, evaluated over covariates, treatment, and $\hat{\pi}$. This variation allows us to isolate the benefits of the additive structure of $k_{A \cdot PRFB}$.

Radial Basis Function (RBF): This is exactly $k_{RBF}$ as described in 3.2, evaluated over covariates, treatment, and $\hat{\pi}$. This variation allows us to isolate the effect of incorporating uncertainty in $\pi$ into a non-additive kernel via $k_{\hat{\pi}}^{PRBF}$.

RBF without dosage (RBF-ND): The same as $k_{RBF}$, but excluding $\hat{\pi}$ entirely, allowing us to study the benefit of conditioning on them.

Additive Symmetrized Gaussian KL-Divergence (A-SymG) [21]: Moreno et al. propose a kernel based on the symmetrized Kullback-Leibler (KL) divergence between Gaussian-distributed inputs (which we call $k_{\hat{\pi}}^{SG}$) to account for uncertainty in inputs and report improvement in classification tasks using support vector machines. See Appendix B for $k_{\hat{\pi}}^{SG}$'s definition. $k_{A \cdot SymG}$ is the same as $k_{A \cdot PRFB}$, but with $k_{\hat{\pi}}^{PRBF}$ replaced with $k_{\hat{\pi}}^{SG}$ over covariates and PS estimates for comparability.

### 4.3 Non-GP models

We compare against popular methods for estimating ADRFs, available in the R library "causaldrf". We estimate ADRFs at observed dosages rather than, say, a grid of evenly-spaced points, allowing us to calculate comparable performance statistics. We provide a full list of methods tested and their specifications in Appendix C. Here, summarize only the two top-performing methods:

Hirano-Imbens Imputation (HI) [17]: OLS regression is used to model $\hat{\pi}$, which in turn are used to fit to a linear function of quadratic treatment, covariate, and dosage terms. As in [16], plug various values for *t* into the fit model, then average responses to estimate the ADRF.

Nonparametric Partial Mean (NPM) [7]: This method uses a GLM to estimate the propensity function. A normal product kernel over treatments and dosage is used to estimate $E[Y|T = t, \Pi = \pi]$; these estimates are averaged over observations to estimate the ADRF.

### 4.4 Full simulation

We adapt the simulations presented by [15] to simulate a continuous treatment. Due to space constraints, Table 1 only provides results from the most difficult simulations: small sample size with a non-linear $\mu$. We include all GP model results and only the top-performing non-GP models. See Appendix D.2 for data-generating processes and the full table of experimental results.

**Table 1: Full Simulation Results – Nonlinear $\mu$**

| | | Homogeneous effect | | | | Heterogeneous effect | | | |
|---|---|---|---|---|---|---|---|---|---|
| *n* | Method | $Cov_{90}$ | $I_{90}$ | Bias | RMSE | $Cov_{90}$ | $I_{90}$ | Bias | RMSE |
| | A-PRBF | **0.95** | 3.21 | 0.181 | 1.18 | **0.91** | 6.30 | -0.177 | 2.79 |
| | A-RBF | 0.95 | 3.21 | 0.181 | 1.18 | 0.91 | 6.32 | -0.179 | 2.78 |
| | PRBF | 0.95 | 3.23 | 0.182 | 1.17 | 0.91 | 6.38 | -0.226 | 2.71 |
| 250 | RBF | 0.77 | 3.28 | 0.175 | 1.18 | 0.87 | 6.57 | 0.235 | 2.46 |
| | RBF-NT | 0.77 | 3.34 | 0.179 | 1.19 | 0.86 | 6.55 | 0.234 | 2.47 |
| | HI | 0.48 | 2.89 | -0.014 | 1.43 | 0.42 | 2.67 | -0.763 | 3.10 |
| | NPM | 0.51 | 3.63 | 0.169 | 1.48 | 0.54 | 3.91 | -0.508 | 2.73 |

## 4.5 IHDP simulation

We simulate a response variable using a real continuous treatment and real covariates from the Infant Development & Health Program study [14]. Specifically, we use the number of days in the Child Development Center as a continuous treatment. We focus specifically on children with non-zero CDC days. We test configurations with homogeneous and heterogenous treatment effects for each of the three simulated response surfaces described in [16]. Table 2 includes results only from the most difficult simulations (surface C) and top-performing non-GP models. See Appendix D.3 for generating processes and the full table of experimental results.

**Table 2: IHDP simulation results – Response surface C**

| Surface | Method | Homogeneous effect | | | | Heterogeneous effect | | | |
|---|---|---|---|---|---|---|---|---|---|
| | | $Cov_{90}$ | $I_{90}$ | Bias | RMSE | $Cov_{90}$ | $I_{90}$ | Bias | RMSE |
| C | A-PRBF | **0.60** | 62.37 | -12.08 | 54.09 | **0.54** | 52.65 | -19.95 | 70.58 |
| | A-RBF | 0.60 | 62.80 | -11.66 | 53.57 | 0.52 | 51.18 | -19.56 | 69.51 |
| | PRBF | 0.57 | 62.64 | -14.56 | 62.30 | 0.53 | 55.86 | -14.13 | 82.71 |
| | RBF | 0.43 | 51.22 | -7.69 | 52.15 | 0.53 | 39.05 | -6.429 | 53.13 |
| | RBF-NT | 0.43 | 51.31 | -7.65 | 53.14 | 0.53 | 39.16 | -6.431 | 53.12 |
| | A-SymG | 0.35 | 49.88 | -12.06 | 55.31 | 0.23 | 48.62 | -11.41 | 54.58 |
| | HI | 0.44 | 69.54 | -22.80 | 48.62 | 0.27 | 71.54 | -15.45 | 61.51 |
| | NPM | 0.69 | **104.8** | -17.53 | 46.43 | 0.38 | 112.1 | -5.880 | 58.79 |

## 4.6 Empirical findings

Findings were consistent between the studies in 4.4 and 4.5, so we summarize observations from both here. See Appendix D for further, problem-specific discussions.

Beginning with ablation findings, we find evidence of RIC and that our priors mitigate it. Including $\hat{\pi}$ in RBF provides no advantage over RBF-NT, whereas including it in PRBF and A-RBF provided a significant advantage over RBF-NT. Whereas RBF does not regularize $\hat{\pi}$ and $t$'s respective kernel parameters independently (GPyTorch does not yet support this), PRBF and A-RBF do (albeit by different mechanisms, described in 3.3 and Appendix D) and can more readily discover heterogeneous dose-response effects. We find that A-PRBF and A-RBF tend to achieve comparable $Cov_{90}$ scores, but A-PRBF yields smaller $I_{90}$s *without* substantive increases in RMSE or Bias.

RBF and RBF-NT models dramatically underestimated interval lengths in all Linear experiments (Table D.1), even after regenerating the data several times. $Y$ itself varies less in these experiments, so we hypothesize that the regularized marginal likelihood can be maximized by collapsing the posterior to vary closely around the ADRF (hence these models still have small bias and RMSE). In nonlinear experiments, mean interval lengths, bias, and RMSE are comparable across GP models, but the average coverage RBF and RBF-NT is much lower. Note that strong confounding occludes the contribution of $t$ to variation in the response, particularly relative to $\pi$; RBF and RBF-NT models lack mechanisms for controlling and differentiating covariation along $t$ and covariation along the $\pi$ dimension. We suspect this problem does not arise for PRBF and A-RBF because they include estimated variation in $\hat{\pi}$, which regularizes covariance with "further" points along $\hat{\pi}$: their lower bound on $\hat{\pi}$'s length-scale controls $\pi$'s contribution to variation in $y$. We find that PRBF performs well in the homogeneous treatment effect experiments relative A-PRBF, but that this difference is reversed in heterogeneous settings. These findings agree with [15]'s description of RIC and argument that independent regularization of $\pi$ and $t$ is necessary for uncovering heterogeneous effects. Regularizing priors over A-PRBF's summed kernel terms prefer larger-normed GP covariance matrix than the PRBF prior would; roughly speaking, we sum prior densities over two scale parameters ($\gamma$ and $\omega$) rather than one ($\gamma$). In our experiments, this assumption is appropriate when we expect greater variation in the response surface (e.g., as a result of more variation in treatment effect), but less so when $\mu$ is linear and/or ATE is homogeneous. This suggests that, in practical situations where less heterogeneity in treatment effects is expected or the ADRF is expected to be a relatively smooth process, we should prefer PRBF by default, or adjust A-PRBF scale parameter priors to reflect that belief.

Overall, we find that GP models with A-PRBF and PRBF kernels provide low-bias, high coverage estimates of the ADRF, relative to all non-GP models tested. We observe that bootstrapped non-GP intervals typically underestimate variation in the ADRF: averaging across experiments, our kernels' GP posteriors (at a given level of $t$) are 30% more likely to include the true ADRF and average $I_{90}$ 131% larger than the best non-GP models. Meanwhile, GP RMSEs are only 9% larger than the best non-GP ADRF estimates; this, taken with the observation that bias results are comparable between our GP kernels and non-GPs methods, indicates that our GPs are not overestimating variation but rather that the non-GP methods' bootstrapped intervals are underestimating it. We should not, however, argue that GP models should be analysts' universal preference—we find non-GP methods, particularly Nonparametric Partial Mean [7], manage competitive coverage rates with much narrower intervals. That said, best non-GP model seemed to depend heavily on the data-generating process, whereas A-PRBF and PRBF models required little tuning across experiments. We believe this constitutes a significant benefit in reducing the risk of misspecification, but nonetheless advise caution to practitioners and recommend a "panel of experts" approach in high-risk scenarios.

Finally, we find A-SymG models consistently underperformed A-PRBF models, reflected in low coverage rates or excessively large interval lengths. This is significant, because, as discussed in in-depth in Appendix A.2, the theory behind $k_{\hat{\pi}}^{PRBF}$ suggests that it allows us to relate its argument distributions to each other *in the context* of GP prior beliefs over response hypotheses. Meanwhile, a kernel based on Kullback-Leibler divergence broadly captures divergence in pointwise posteriors for $\pi$ without regard their weighting of hypothesis functions under the GP kernel.

# 6 Conclusion

We find GP regression to be a versatile, accurate approach for analyzing dose-response relationships in the observational setting. We observe competitive performance with little tuning, while the non-GP method performance varied considerably from one data-generating process to another. Beyond these practical benefits, we identified the mechanism by which RIC occurs in the GP response model and find that naïve inclusion of a propensity score in a Gaussian process response model offers benefits only up to our ability to independently regularize treatment parameters relative to it. In our case, as in Hahn et al. [15], we identify priors that address model-specific forms of RIC, suggesting that other, existing observational causal inference methods may stand to benefit from similar mitigations. Finally, we demonstrate the value in incorporating uncertainty in dosage propensities via $k_{\hat{\pi}}^{PRBF}$, and believe that further generalizations and connections to the kernel embedding literature stand to benefit GP modeling in a broader range of causal inference problems and datasets.

We have applied this approach at Amazon to make better decisions about how to allocate resources, as part of broader analysis efforts which first involve feature selection, then ADRF modeling, then more precise estimation of (conditional) ATEs using approaches like Double ML estimation [4]. More broadly, we find ADRF estimation to be a valuable yet often overlooked tool in the analyst's toolbox, which can further research in cases where experiments are not feasible and units select into varied *levels* of treatment. We hope this contribution marks a step toward drawing researchers' attention to such problems and approaches, allowing for richer descriptions of a complex world.

**Acknowledgement**: I thank Carl-Johann Simon-Gabriel for lending his expertise and time to a thoughtful, thorough review of this work.

# Appendix A – Technical Results

## A.1 – GP regressor as the best linear unbiased predictor

Define the best linear unbiased predictor (BLUP) to be the predictor $\hat{Y} = \sum_{i=1}^{n} c_i k(x_i)$ that minimizes the variance of prediction error f for new observation $Y'$,

$$V = Var(Y' - \hat{Y}) \propto \mathcal{N}(\hat{Y}|Y').$$

Suppose $Y(x)$ is a Gaussian process such that, for every finite set of inputs $x_1, \ldots, x_n$,

$$Y(x)|_{x=x_1,\ldots,x_n} = (Y(x_1), \ldots, Y(x_n)) = y,$$

is a multivariate Gaussian random variable with $k(\cdot,\cdot)$, a real positive-definite kernel defined over $x \in \mathbb{R}^d$, comprising its covariance matrix $K$. By Moore-Aronszanj Theorem, because $k$ is a positive definite kernel, there exists a corresponding Reproducing Kernel Hilbert Space (RKHS), $\mathcal{H}$.

Consider a functional $J[f] = \frac{\alpha}{2} [\![f]\!]_{\mathcal{H}}^2 + Q(y, f)$ over members of $\mathcal{H}$, such that $f \in \mathcal{H}$, Q is a quasi-likelihood (i.e., the negative log likelihood over our $n$ training points up to some terms not involving $f$), and $f = (f(x_1), \ldots, f(x_n))$. The Representer Theorem shows that the minimizer of $J$ has the form $f(x) = \sum_{i=1}^{n} c_i k(x, x_i)$.

[22] show that the Gaussian process regressor predictive mean (for a new point $x'$) also satisfies the form $E[f'|y] = \sum_{i=1}^{n} c_i k(x', x_i)$, and that, when $Q$ is the Gaussian likelihood, the coefficients minimizing $J$ are $(c_1, \ldots, c_n) = (K + \sigma_n^2 I)^{-1} y$. Substituting the coefficients in, we see the minimizing function of $J$, $\hat{f'} = k(x')(K + \sigma_n^2 I)^{-1} y$, both satisfies linearity and minimizes the variance of prediction error when the regressor and modeled process kernels coincide.

## A.2 – Predictive Radial Basis Function kernels

Consider the Gaussian kernel,

$$k_l(x, x') = \frac{1}{l\sqrt{2\pi}} \exp\left(-\frac{\|x - x'\|^2}{2l^2}\right) = \mathcal{N}(x \mid x', l^2),$$

for $x \in \mathcal{X}$. Typically, we define a Gaussian process $\mathcal{GP}$ such that

$$f(x) \sim \mathcal{GP}(m(x), k(x, x')),$$

where $m$ is an arbitrary mean function, and $k$ is defined such that

$$k(x, x') \coloneqq E[(f - m)(f' - m')|x, x'],$$

where we write $f(x) = f$ and $f(x') = f'$ for brevity (likewise for $m$).

We consider instead the case where an input $\tau$ is generated by a Gaussian distribution whose mean and variance depend on $x \in \mathcal{X} \subset \mathbb{R}^d$, where $\mathcal{X}$ is an open subset (we dispense with bolded font for readability, as all input spaces are assumed to be vector spaces); that is,

$$\tau \sim Q(x) = \mathcal{N}(\mu(x), \sigma(x)).$$

We assume $|\mu| < \infty$ and $0 < \sigma < \infty$ for all $x \in \mathcal{X}$. Now we assume, without loss of generality, that $m(\theta) = 0$ (see [23] for handling of non-zero mean cases). We would like to account for uncertainty or known variation in $\pi$. For each pair of observations $x, x' \in \mathcal{X}$, say we choose to incorporate uncertainty by averaging kernels evaluated over respective samples from $Q(x)$ and $Q(x')$, defining

$$k_\pi(x, x') \coloneqq E_{\pi \sim Q(x)}\big[E_{\pi' \sim Q(x')}[E[ff'|\pi, \pi']|x]|x'\big]. \tag{8}$$

In the case that $E[ff'|\pi, \pi']$ is the Gaussian kernel, we can write these expectations as

$$k_\pi^{PRBF}(x, x') = \int\int \frac{1}{\sqrt{2\pi l^2}} \exp\left(-\frac{(\theta - \theta')^2}{2l^2}\right) \frac{1}{\sqrt{2\pi \sigma(x)^2}} \exp\left(-\frac{(\mu(x) - \theta)^2}{2\sigma(x)^2}\right) d\theta \, \frac{1}{\sqrt{2\pi \sigma(x')^2}} \exp\left(-\frac{(\mu(x') - \theta')^2}{2\sigma(x')^2}\right) d\theta'$$

$$= \int \frac{1}{\sqrt{2\pi(l^2+\sigma^2)}} \exp\left(-\frac{(\theta'-\mu)^2}{2(l^2+\sigma^2)}\right) \frac{1}{\sigma'\sqrt{2\pi}} \exp\left(-\frac{(\mu'-\theta')^2}{2\sigma'^2}\right) d\theta'$$

$$= \frac{1}{\sqrt{2\pi(l^2+\sigma(x)^2+\sigma(x\prime)^2)}} \exp\left(-\frac{(\mu(x)-\mu(x\prime))^2}{2(l^2+\sigma(x)^2+\sigma(x\prime)^2)}\right) \tag{9}$$

where $\mu = \mu(x)$, $\mu' = \mu(x')$, $\sigma = \sigma(x)$, and $\sigma' = \sigma(x')$.

The expression in (9) is consistent with the result in [12], which Girard derived in the interest recovering covariance conditional only on the input (rather than the input and $\pi$). It is not given that this function is a positive definite (p.d.) kernel (and actually describes a valid GP regressor), nor is it obvious why this method of accounting for uncertainty is useful.

We proceed first by showing that $k_\pi$ is indeed a valid p.d. kernel, then that it is particularly useful because it corresponds to inner products between distribution embeddings into a Reproducing Kernel Hilbert Space (RKHS) [13, 29, 30]. Using this fact, we show that the p.d. matrix $K_\pi \subset \mathbb{R}^{n \times n}$ with entries $k_\pi(x_i, x_j), i, j \in \{1, \dots, n\}, n \in \mathbb{N} \setminus \{0\}$ then contains all necessary information to compute Maximum Mean Discrepancy.

### A.2.1 Proof of positive-definiteness

In regards to identification as a p.d. kernel, $k_\pi$ appears similar to the Gaussian kernel but differs fundamentally in structure: in particular, $\sigma$ and $\mu$ vary as functions of its arguments, and cannot be decomposed into univariate functions using techniques previously applied in the literature [35]. Instead, we focus on translation-invariant kernels, identify them as positive-definite functions, then prove that inner products of translation-invariant kernel mean embeddings are themselves positive-definite kernels, proving $k_\pi$'s positive-definiteness.

**Preliminary** Let $\Omega \subset \mathbb{R}^d$ denote an open set and $F: \Omega - \Omega \to \mathbb{C}$ a function, where $\Omega - \Omega = \{x - y : x, y \in \Omega\}$. $F$ is a positive definite function if, for any $m \in \mathbb{N}, x_1, \dots, x_m \in \Omega$, and $c_1, \dots, c_m \in \mathbb{C}$,

$$\sum_{j,k=1}^{m} F(x_j - x_k) c_j \overline{c_k} \geq 0.$$

**Lemma 1** *Let $(M, \mathfrak{m}, \zeta)$ be a measure space, and let $K: \Omega \times \Omega \to \mathbb{C}$ denote a translation-invariant positive-definite kernel such that $K(x_j, x_k) = K(x_j - x_k, 0)$. Then, $K$ is a positive-definite function.*

**Proof** Using basic properties of Reproducing Kernel Hilbert Spaces, we note any positive definite kernel can be written as

$$K(x_j, x_k) = \int_M k(x_j, s) \overline{k(x_k, s)} d\zeta(s),$$

for a $\zeta$-square integrable function $k(x, \cdot): M \to \mathbb{C}$. Letting $x_1, \dots, x_m \in \Omega$, $c_1, \dots, c_m \in \mathbb{C}$, and $F(x_j - x_k) = K(x_j - x_k, 0)$, we have

$$\sum_{j,k=1}^{m} F(x_j - x_k) c_j \overline{c_k} = \sum_{j,k=1}^{m} \int_M k(x_j, s) \overline{k(x_k, s)} d\zeta(s)$$

$$= \int_M \left| \sum_{j=1}^{m} k(x_j, s) c_j \right|^2 d\zeta(s) \geq 0,$$

or rather, that $K$ is a positive definite function. ∎

**Preliminary** For a distribution $D \in C_0^\infty(\Omega)$, define the kernel mean embedding $F_D: \Omega \to \Omega$ under a positive-definite kernel $K$ as $F_D(\theta) = \int K(\cdot, \theta) \, dD(\theta)$. The inner product of kernel mean embeddings is then

$$\langle F_D, F_T \rangle = \langle \int K(\cdot, \theta) \, dD(\theta), \int K(\cdot, \theta') \, dT(\theta') \rangle = \int \int K(\theta, \theta') dD(\theta) dT(\theta').$$

**Theorem 1** *Let $K: \Omega \times \Omega \to \mathbb{C}$ be a translation-invariant positive definite kernel (and thus a positive-definite function, per Lemma 1), and let $\varphi, \psi \in C_0^\infty(\Omega)$. Then, $\langle F_\varphi, F_\psi \rangle$ itself is a positive-definite kernel.*

**Proof** Recall that, due to $K$'s translation invariance, we can write the inner product of kernel mean embeddings as

$$\langle F_\varphi, F_\psi \rangle = \int_\Omega \int_\Omega F(\theta - \theta') \, \varphi(\theta) \overline{\psi(\theta')} d\theta d\theta'.$$

Then, we have, for any $\varphi, \psi \in C_0(\mathbb{R})$,

$$\sum_{j,k=1}^{m} c_j c_k \langle F_\varphi, F_\psi \rangle = \sum_{j,k=1}^{m} c_j c_k \int_\Omega \int_\Omega F(\theta - \theta') \, \varphi(\theta) \overline{\psi(\theta')} d\theta d\theta'$$

$$= \int_\Omega \int_\Omega \sum_{j,k=1}^{m} c_j \varphi(\theta) c_k \overline{\psi(\theta')} \, F(\theta - \theta') \, d\theta d\theta' = \int_\Omega \left| \sum_{j=1}^{m} F(\theta - \theta') a_j(\theta) \right|^2 d\theta \geq 0$$

for all $c_1, \ldots, c_m \in \mathbb{R}$ and where $a_j(\theta) = c_j \varphi(\theta)$, proving that $\langle F_\varphi, F_\psi \rangle$ is itself a positive-definite kernel. ∎

By recognizing the Gaussian kernel as translation-invariant and identifying (9) with $\langle F_\varphi, F_\psi \rangle$, we find that $k_\pi$ is subject to Theorem 1 and is a positive-definite kernel.

### A.2.2 Connection to kernel mean embeddings

To understand why this particular treatment of uncertainty is useful, we connect the kernel in (9) with the literature on kernel mean embeddings (KMEs) [13, 29, 30]. We find that $k_\pi(x, x')$ is a closed form expression for the inner product between RKHS embeddings of distributions $\varphi_x$ and $\psi_{x'}$. $k_\pi(x, x')$ achieves a lower bound in the case that $\varphi_x$ and $\psi_{x'}$ are identical, which coincides with the case that the so-called Maximum Mean Discrepancy (MMD) is zero. Intuitively, MMD upper-bounds the difference in mean values of all functions in the RKHS when means are taken with respect to either $\varphi_x$ and $\psi_{x'}$; we will make this precise after a few preliminaries.

Let $D$ denote a probability measure and $k$ an arbitrary p.d. kernel (both of which take values over $\Omega$, an open subset of $\mathbb{R}^d$). In this case, $\int_\Omega k(\cdot, \theta) \, D(\theta) d\theta = \int_\Omega k(\cdot, \theta) \, dD(\theta) = \phi_D$ is called the KME of $D$. Simon-Gabriel and Schölkopf [29] generalize this to Schwartz-distributions, the set of continuous linear functionals on $C_c^\infty(\Omega)$, denoted $\mathcal{D}^\infty$, which contains all continuous functions and signed measures. The authors prove that, for distributions $D, T \in \mathcal{D}^\infty$,

$$\langle \phi_D, \phi_T \rangle = \langle \int k(\cdot, \theta) \, dD(\theta), \int k(\cdot, \theta') \, dT(\theta') \rangle = \int \int k(\theta, \theta') dD(\theta) dT(\theta') \tag{10}$$

Clearly, as defined in (10), $\varphi_x(\theta), \psi_{x'}(\theta) \in \mathcal{D}^\infty$, yielding the insight that $k_\pi$ is a special case of the inner product between distribution embeddings into $\mathcal{H}_k$.

Maintaining notation, MMD is defined as

$$MMD(D, T, \mathcal{H}_k) \coloneqq \|\phi_D - \phi_T\|_{\mathcal{H}_k} = \sup_{\substack{f \in \mathcal{H}_k: \\ \|f\|_{\mathcal{H}_k} \leq 1}} |Df - Tf|.$$

It is well-known that being able to compute (10) allows computation of the MMD; Gretton et al. [13] show by way of the reproducing property that

$$MMD^2(D, T, \mathcal{H}_k) = E_{\theta, \theta' \sim D}[k(\theta, \theta')] - 2E_{\theta \sim D, \delta \sim T}[k(\theta, \delta)] + E_{\delta, \delta' \sim T}[k(\delta, \delta')]$$

$$= \langle \phi_D, \phi_D \rangle - 2\langle \phi_D, \phi_T \rangle + \langle \phi_T, \phi_T \rangle,$$

where all samples from $D$ and $T$ are independent. Identifying inner products with $k_\pi$, we see

$$k_\pi(x,x) = \langle \phi_{\varphi_x}, \phi_{\varphi_x} \rangle = \|\varphi_x\|_{\mathcal{H}_k} = \frac{1}{\sqrt{2\pi(l^2 + 2\sigma(x)^2)}},$$

so we write

$$MMD^2(\varphi_x, \psi_{x\prime}, \mathcal{H}_k) = \|\varphi_x\|_{\mathcal{H}_k} + \|\psi_{x\prime}\|_{\mathcal{H}_k} - 2k_\pi(x, x'),$$

Per Gretton [13], Simon-Gabriel and Schölkopf [29], $MMD^2(\varphi_x, \psi_{x\prime}, \mathcal{H}_k) = MMD(\varphi_x, \psi_{x\prime}, \mathcal{H}_k) = 0$ iff $\varphi_x = \psi_{x\prime}$, which coincides with the case where $k_\pi(x, x') = k_\pi(x, x)$. Note that a kernel matrix consisting of $k_\pi(x_i, x_j), i, j \in \{1, ..., n\}, n \in \mathbb{N} \setminus \{0\}$ then contains all necessary information for the MMD between any particular pair of embeddings.

Ritter [25] offers an additional perspective by showing that, in the case where functions are sampled from a zero-mean GP,

$$MMD^2(\varphi_x, \psi_{x\prime}, \mathcal{H}_k) = E_{f \sim \mathcal{GP}(0,k)}[(\varphi_x f - \psi_{x'} f)^2] = var[\varphi_x f - \psi_{x'} f],$$

which implies $(\varphi_x f - \psi_{x'} f) \sim \mathcal{N}(0, MMD(\varphi_x, \psi_{x'}, \mathcal{H}_k))$. That is, up to $var[\varphi_x f - \psi_{x'} f] - \|\varphi_x\|_{\mathcal{H}_k} - \|\psi_{x'}\|_{\mathcal{H}_k}$, $k_\pi$ allows us to express covariance as the concurrence in RKHS embeddings of a GP prior's hypotheses under each pair of pointwise distributions $\varphi_x, \psi_{x'}$. This is a powerful intuition about $k_\pi$'s role in GP regression: we can easily relate distributions to each other *in the context* of GP prior beliefs over hypotheses; this is in contrast to alternatives such as Kullback-Leibler divergence (as in $k_{SymG}$, see Appendix B) [21], where the relationship between kernel, members of the RKHS, and GP prior beliefs are not clear.

### A.3 – Double Robustness of the Additive Estimator

The flexibility of Gaussian Process modeling and the additive construction proposed in this work allow us to prove convergence under either a correctly specified response function or a correctly specified generalized propensity score. We sketch such a proof for the model specification given in (7); to satisfy double robustness, we require

$$E\left[f_{true}\left((\pi_{true}(X), X)\right) + g_{true}(T) - f\left((\hat{\pi}(X), X)\right) + g(T)\right] \to 0$$

when either $f$ and $g$ are correctly specified or $\hat{\pi}$ is.

First, we consider the case where the response function is specified correctly but the generalized propensity score is not (i.e., $\pi \neq \hat{\pi}$); then, $E[g(T)] \to g_{true}(T)$ and $E[f((\hat{\pi}(X), X))] \to E[f_{true}((\hat{\pi}(X), X))]$ due to the consistency of GP regression. We then require only that $E[f_{true}((\hat{\pi}(X), X))] = E[f_{true}((\pi_{true}(X), X))]$, which is satisfied asymptotically if $\hat{\pi} \perp T \mid X$, a given since $\hat{\pi}$ does not contain any information about $T$ not already in $X$; we only include both for efficiency.

Now, for the case that $\hat{\pi} = \pi_{true}$ is correctly specified but $f_{true} \neq f$ and $g_{true} \neq g$, we require

$$Bias(t) = E\left[f\left((\hat{\pi}(X), X)\right) + g(t) \mid T = t, \hat{\pi} = \pi_{true}\right] - Y(t) \to 0.$$

it is the case that $P(T = t \mid X, \hat{\pi}(X)) = P(T = t \mid \hat{\pi}(X))$, the so-called balancing property. This property implies $E[f((\hat{\pi}(X), X)) \mid T = t] = E[f_{true}((\hat{\pi}(X), X))]$. Applying this, and the law of iterated expectations and GP regression consistency so that $E[g(t)] = g_{true}(t)$, we have

$$Bias(t) = E\left[f_{true}\left((\hat{\pi}(X), X)\right)\right] + g_{true}(t) - Y(t).$$

We assume that $f_{true}\left((\hat{\pi}(X), X)\right)$ captures all confounding, so that

$$E\left[f_{true}\left((\hat{\pi}(X), X)\right)\right] = Y(t) - g_{true}(t),$$

or

$$Bias(t) = g_{true}(t) - \left(Y(t) - g_{true}(t)\right) - Y(t) = 0.$$

This sketch is given purely for illustration but similar logic can be applied for other model specifications, such as $f\big((\hat{\pi}(X), X)\big) = f_1\big(\hat{\pi}(X)\big) + f_2(X)$.

## Appendix B – A-SymG kernel formulation and comparison

The Symmetrized Gaussian KL-Divergence (SymG) kernel was originally developed for Support Vector Classification over feature spaces with Gaussian-distributed noise. We study it relative

As a first step, [21] a Gaussian distribution $p(x|\theta_i)$ for each observation $i$, then use $p(x|\theta_i)$ in

$$K_{SG}\big(p(x|z_i), p(x|z_j)\big) = exp\left(-\frac{D_{SG}\big(p(x|z_i), p(x|z_j)\big)}{A}\right),$$

$$D_{SG}\big(p(x|z_i), p(x|z_j)\big) = \int_{-\infty}^{\infty} p(x|z_i) \log\left(\frac{p(x|z_i)}{p(x|z_j)}\right) dx + \int_{-\infty}^{\infty} p(x|z_j) \log\left(\frac{p(x|z_j)}{p(x|z_i)}\right) dx$$

$$= tr(\Sigma_i \Sigma_j^{-1}) + tr(\Sigma_j \Sigma_i^{-1}) + tr\left((\Sigma_i^{-1} + \Sigma_j^{-1})(\mu_i - \mu_j)(\mu_i - \mu_j)^T\right) - 2d,$$

where $z_i = (\mu_i, \Sigma_i)$, the tuple of the estimated mean and covariance matrix, respectively, for observation $i$, $d$ is the dimension of $x$, and A is a parameter. In our experiments, we use

$$k_{A \cdot SymG}(\theta, \theta') = k_{\hat{\pi}}^{SG}\big((\hat{\pi}(x), x), (\hat{\pi}(x'), x')\big) + k_{RBF}(t, t'),$$

allowing us to directly compare the $k_{\hat{\pi}}^{SG}$ mechanism for incorporating uncertainty to $k_{\hat{\pi}}^{PRBF}$. Due to time constraints, we only report results for SymG in the IHDP simulation (Appendix D.2).

Experimental results indicate that $k_{A \cdot SymG}$ significantly underperforms $k_{A \cdot PRBF}$ in terms of coverage and interval length (see Appendix D.2), but we have yet to understand exactly why. In a future version of this work, we will look to study the positive-definiteness, RKHS, and behavior of this kernel in more depth.

## Appendix C – Non-GP model specifications

Here we provide the full list of non-GP model definitions and provide detail on hyperparameters used, where applicable:

Additive Splines (ADS) [3]: This method estimates generalized propensity scores using generalized linear models (GLM). Scores are then input to an additive spline model including pairwise interactions of treatment and propensity spline terms for each user-defined knot. We use 3 knots as this seemed to provide the best tradeoff between model fit and smoothness.

Bayesian additive regression trees (BART) [5, 16]: Iterative fitting of piecewise-constant regression trees, each comprised of decision nodes partitioning the covariate space. Trees are summed to form a single response surface. [5] propose a priors over tree depth and the fit of tree nodes. While not explicitly specified, [16] imply that the ADRF can be found by plugging in continuous values of $t$ to the updated model and averaging. We do not report $Cov_{90}$ or $I_{90}$ for BART models as their "causaldrf" implementation does not expose posteriors and one update takes approximately as long as the bootstrap procedure of non-GP methods, so bootstrapping was not feasible. Using the notation in [16], we set $k$=2, use 3 degrees of freedom in the chi-squared prior over $\sigma$, set $q = 0.9$, and the default tree prior choices described in [5] (as these are found to be immaterial to fit [16]). For all BART models, we draw 1100 Markov Chain Monte Carlo samples; we do not bootstrap performance statistics as one update takes approximately as long as the bootstrap procedure of non-GP methods.

BART with GBM Propensity Scores (BART-PS-GBM) [15]: Incorporates an estimate of propensity score into BART to protect against RIC. [14] Report a dramatic reduction in RIC but do not specify a preferred propensity score estimator but we model it using a Gradient Boosting Machine with up to thirty trees (with the actual number determined by cross-validation over the propensity score holdout set), a shrinkage parameter of 0.1, and an interaction depth of 3. [15] report dramatic reduction of RIC. We use the same prior settings described for BART. We do not use the Bayesian Causal Forest formulation or priors due to their specification to the binary treatment setting.

BART with OLS Propensity Scores (BART-PS-OLS) [15]: The same as BART-PS-GBM but OLS linear regression is used to model propensity scores.

Hirano-Imbens Imputation (HI) [17]: OLS regression is used to model expected dosage, which in turn are used to fit to a linear function of quadratic treatment, covariate, and dosage terms. As in [16], plug various values of $t$ into the fit model, then average responses to estimate the ADRF.

Nadaraya-Watson (NW) [7]: This method is a kernel smoothing of an inverse propensity weighted estimator for $y$ at a given $t$. In this case the propensity is estimated a GLM extension of the propensity model used in [17].

Nonparametric Partial Mean (NPM) [7]: This method uses a GLM to estimate the propensity function. A normal product kernel over treatments and propensity scores is used to estimate $E[Y|T=t, \Pi=\pi]$; these estimates are averaged over observations to estimate the ADRF.

Inverse Propensity of Treatment Weighting (IPTW) [26]: This is a continuous adaptation of the Horvitz-Thompson-type estimators. GLM regression is used to the estimate treatment likelihood function parameters; the ratio of the treatment likelihood and the conditional likelihood given the parameters is used to reweight the regression of $Y$ on $t$.

## Appendix D – Simulation definition, results, and analysis

### D.1 – Evaluation metrics

In all studies, we sample estimates of the bias of the outcome model relative to $E[Y(t)|T=t] = \tau(t)$,

$$\widehat{Bias}_j = \frac{1}{n}\sum_{i=1}^{n} \tau(t_i) - \hat{Y}_j(\boldsymbol{\theta}_i),$$

by sampling $\hat{Y}$ from our outcome model (here $j$ indexes $p$ samples drawn from the outcome model posterior and $\tau(t_i)$ is the true DRF at the level of treatment observed for that unit). We also report an estimate of the proportion of $\tau(t)$'s covered by a 90% credible interval over the GP posterior,

$$\widehat{Cov}_{90} = \frac{1}{n}\sum_{i=1}^{n} \mathbb{I}\left(\hat{q}_5(\boldsymbol{\theta}_i) \leq \tau(t_i) \leq \hat{q}_{95}(\boldsymbol{\theta}_i)\right), \text{ and}$$

$$\widehat{I_{90}} = \sum_{i=1}^{n} \widehat{q_{95}}(\boldsymbol{\theta}_i) - \widehat{q_5}(\boldsymbol{\theta}_i),$$

where $\mathbb{I}$ is the indicator function and $\widehat{q_r}$ is the r-percentile of the Gaussian outcome GP posterior. Note that $\hat{\mu}(\boldsymbol{z}_i)$ and $\hat{\sigma}(\boldsymbol{z}_i)$ are deterministic functions so we do not need to sample from the posterior to estimate $\widehat{Cov}_{90}$. $\widehat{I_{90}}$, the mean credible interval length, allows us contextualize $\widehat{Cov}_{90}$ estimates.

In addition to these metrics, we also estimate the root mean squared error of the model with respect to the true means underlying the observed outcomes,

$$\widehat{RMSE}_j = \sqrt{\frac{1}{n}\sum_{i=1}^{n}\left(\tau(t_i) - \hat{Y}_j(\boldsymbol{\theta}_i)\right)^2},$$

again sampling from the GP posterior. We report the mean of $\widehat{Bias}$ and $\widehat{RMSE}$, taken over all $j \in \{1, \ldots, p\}$ samples.

For non-GP methods (except BART, as explained below), we similarly estimate all four of these statistics by bootstrap sampling 50 training sets.

### D.2 – Full Simulation

We adapt the simulation presented by [15] to simulate a continuous treatment:

$$y_i = \mu(\pmb{x}_i) + \eta_i \varphi(\pmb{x}_i), \text{ where}$$

$$\mu(\pmb{X}) = \begin{cases} 1 + g(X_5) + X_1 X_3, & linear \\ 1 + g(X_5) + 6|X_3 - 1|, & nonlinear \end{cases},$$

$$\eta_i \sim \mathcal{N}(\pi(\pmb{x}_i), 1),$$

$$\pi(\pmb{x}_i) = 0.8\Phi\left(\frac{3\mu(\pmb{x}_i)}{s} - \frac{x_{i,1}}{2}\right) + \frac{u_i}{10}, \quad u_i \sim Uniform(0,1),$$

$$\varphi(\pmb{X}) = \begin{cases} 3, & homogeneous \\ 1 + 2X_2 X_4, & heterogeneous \end{cases},$$

$$g(x) = \begin{cases} 2, & x = 1 \\ -1, & x = 2 \\ -4, & x = 3 \end{cases},$$

$X_{1,2,3} \sim \mathcal{N}(0,1)$, $X_4 \sim Bernoulli(0.5)$, $X_5 \sim Multinoulli((0.1, 0.4, 0.5))$, and $s$ is the standard deviation of $\mu$ taken over the simulated sample. We alternate options for $\mu(\pmb{X})$ and $\varphi(\pmb{X})$, and alternate between samples of 250 and 500 observations for a total of eight distinct simulations. We provide results in Tables D.1 and D.2:

Table D.1: Full Simulation Results – Linear $\mu$

| $n$ | Method | Homogeneous effect | | | | Heterogeneous effect | | | |
|---|---|---|---|---|---|---|---|---|---|
| | | $Cov_{90}$ | $I_{90}$ | Bias | RMSE | $Cov_{90}$ | $I_{90}$ | Bias | RMSE |
| 250 | A-PRBF | 0.90 | 3.57 | -0.097 | 1.05 | **1.00** | 7.42 | -0.204 | 1.15 |
| | A-RBF | 0.89 | 3.56 | -0.096 | 1.04 | 1.00 | 7.78 | -0.215 | 1.12 |
| | PRBF | **0.91** | 3.57 | -0.078 | 0.94 | 1.00 | 7.58 | -0.193 | 1.03 |
| | RBF | 0.12 | 0.22 | -0.083 | 0.69 | 0.50 | 0.89 | -0.135 | 0.76 |
| | RBF-NT | 0.13 | 0.24 | -0.081 | 0.69 | 0.46 | 0.86 | -0.136 | 0.77 |
| | ADS | 0.54 | 1.48 | -0.325 | 1.12 | 0.38 | 1.59 | -0.271 | 1.44 |
| | BART | - | - | -0.510 | 1.49 | - | - | -0.176 | 1.24 |
| | BART-PS-GBM | - | - | -0.494 | 1.38 | - | - | -0.270 | 1.38 |
| | BART-PS-OLS | - | - | -0.546 | 1.49 | - | - | -0.168 | 1.27 |
| | HI | 0.42 | 1.31 | -0.368 | 1.23 | 0.41 | 1.55 | -0.423 | 1.36 |
| | NW | 0.05 | 1.42 | -1.034 | 2.90 | 0.26 | 1.30 | -0.380 | 2.15 |
| | NPM | 0.57 | 1.65 | -0.374 | 1.14 | 0.61 | 2.10 | -0.242 | 1.44 |
| | IPTW | 0.17 | 0.90 | -0.479 | 1.33 | 0.38 | 1.26 | -0.195 | 1.38 |
| 500 | A-PRBF | 0.92 | 3.48 | 0.19 | 1.32 | **1.00** | 7.18 | 0.067 | 0.70 |
| | A-RBF | 0.92 | 3.47 | 0.188 | 1.32 | 1.00 | 7.49 | 0.077 | 0.70 |
| | PRBF | **0.93** | 3.42 | 0.171 | 1.26 | 1.00 | 7.44 | 0.119 | 0.72 |
| | RBF | 0.09 | 0.06 | 0.139 | 1.15 | 0.13 | 0.26 | -0.030 | 0.82 |
| | RBF-NT | 0.09 | 0.06 | 0.140 | 1.16 | 0.19 | 0.27 | -0.032 | 0.83 |
| | ADS | 0.62 | 0.97 | -0.049 | 0.50 | 0.70 | 1.34 | 0.109 | 0.79 |
| | BART | - | - | -0.211 | 0.62 | - | - | 0.048 | 0.78 |
| | BART-PS-GBM | - | - | -0.184 | 0.54 | - | - | 0.033 | 0.76 |
| | BART-PS-OLS | - | - | -0.188 | 0.58 | - | - | 0.032 | 0.77 |
| | HI | 0.50 | 0.89 | -0.025 | 0.66 | 0.54 | 1.28 | 0.134 | 0.96 |
| | NW | 0.12 | 0.74 | -0.607 | 1.71 | 0.09 | 0.74 | **-0.030** | 1.44 |
| | NPM | 0.62 | 1.09 | -0.086 | 0.48 | 0.77 | 1.43 | 0.169 | 0.81 |
| | IPTW | 0.40 | 0.63 | -0.203 | 0.54 | 0.38 | 1.00 | 0.092 | 0.90 |

When $\mu$ is linear, we observe comparable results between A-PRBF, A-RBF, and PRBF. As we would expect, when treatment effects are homogeneous, there is almost no difference in performance between A-PRBF and A-RBF. However, when treatment effects are heterogeneous, A-PRBF achieves 100% ADRF coverage with a smaller average interval and less bias than A-RBF, suggesting the inclusion of PS uncertainty is useful in right-sizing the posterior variances.

Table D.2: Full Simulation Results – Nonlinear $\mu$

| n | Method | Homogeneous effect | | | | Heterogeneous effect | | | |
|---|---|---|---|---|---|---|---|---|---|
| | | $Cov_{90}$ | $I_{90}$ | Bias | RMSE | $Cov_{90}$ | $I_{90}$ | Bias | RMSE |
| 250 | A-PRBF | **0.95** | **3.21** | 0.181 | 1.18 | **0.91** | 6.30 | **-0.177** | 2.79 |
| | A-RBF | **0.95** | **3.21** | 0.181 | 1.18 | 0.91 | 6.32 | -0.179 | 2.78 |
| | PRBF | 0.95 | 3.23 | 0.182 | **1.17** | 0.91 | 6.38 | -0.226 | 2.71 |
| | RBF | 0.77 | 3.28 | 0.175 | 1.18 | 0.87 | 6.57 | 0.235 | **2.46** |
| | RBF-NT | 0.77 | 3.34 | 0.179 | 1.19 | 0.86 | 6.55 | 0.234 | 2.47 |
| | ADS | 0.51 | 3.22 | -0.076 | 1.47 | 0.42 | 3.00 | -0.591 | 2.99 |
| | BART | - | - | 0.196 | 1.41 | - | - | -0.458 | 3.08 |
| | BART-PS-GBM | - | - | 0.200 | 1.43 | - | - | -0.428 | 3.09 |
| | BART-PS-OLS | - | - | 0.186 | 1.45 | - | - | -0.488 | 3.06 |
| | HI | 0.48 | 2.89 | **-0.014** | 1.43 | 0.42 | 2.67 | -0.763 | 3.10 |
| | NW | 0.22 | 2.04 | 0.361 | 1.89 | 0.24 | 1.83 | -0.768 | 3.67 |
| | NPM | 0.51 | 3.63 | 0.169 | 1.48 | 0.54 | 3.91 | -0.508 | 2.73 |
| | IPTW | 0.37 | 1.94 | 0.148 | 1.46 | 0.33 | 2.62 | -0.733 | 3.44 |
| 500 | A-PRBF | **0.93** | 2.42 | -0.520 | 1.58 | **0.84** | 5.65 | -0.387 | 2.11 |
| | A-RBF | 0.92 | 2.40 | -0.519 | 1.58 | 0.84 | 5.73 | -0.394 | 2.15 |
| | PRBF | **0.93** | 2.39 | -0.519 | 1.58 | 0.81 | 5.64 | -0.394 | 2.10 |
| | RBF | 0.63 | 2.00 | -0.549 | 1.58 | 0.71 | 3.68 | -0.320 | **1.89** |
| | RBF-NT | 0.70 | 2.67 | -0.557 | 1.59 | 0.71 | 3.69 | -0.321 | **1.89** |
| | ADS | 0.70 | 2.35 | -0.137 | 0.96 | 0.50 | 2.31 | 0.011 | 2.08 |
| | BART | - | - | -0.307 | 0.98 | - | - | 0.061 | 2.24 |
| | BART-PS-GBM | - | - | -0.343 | 1.05 | - | - | 0.072 | 2.20 |
| | BART-PS-OLS | - | - | -0.299 | 0.97 | - | - | 0.058 | 2.25 |
| | HI | 0.80 | 2.34 | 0.001 | 1.00 | 0.21 | 2.08 | -0.306 | 2.35 |
| | NW | 0.13 | 1.74 | -0.97 | 2.54 | 0.25 | 1.53 | -0.664 | 3.14 |
| | NPM | 0.80 | 2.63 | -0.04 | 0.86 | 0.53 | 2.77 | 0.003 | 2.10 |
| | IPTW | 0.40 | 1.44 | -0.35 | 0.98 | 0.37 | 1.64 | -0.216 | 2.20 |

In all linear experiments (Table D.1), RBF and RBF-NT models dramatically underestimated interval lengths, even after regenerating the data several times. Y itself varies much less in these experiments, so the regularized marginal likelihood can be maximized by collapsing the posterior to vary closely around the ADRF (hence these models still have small bias and RMSE). In nonlinear experiments, mean interval lengths, bias, and RMSE are comparable across GP models, but the average coverage of RBF and RBF-NT is much lower, suggesting pointwise posterior interval lengths fail to represent variation in the DRF. Note that strong confounding occludes the contribution of $t$ to variation in $Y$, particularly relative to the PS; RBF and RBF-NT models lack mechanisms for controlling and differentiating covariation along $t$ and covariation along the PS dimension. We suspect this problem does not arise for PRBF and A-PRBF because they include estimated variation in $\hat{\pi}$, which regularizes covariance with "further" points along $\hat{\pi}$: their lower bound on $\hat{\pi}$'s length-scale controls the PS's contribution to variation in $Y$. Meanwhile, A-RBF and A-PRBF have explicit priors controlling confounders' contribution to the covariance. In both A-RBF and A-PRBF, priors over the summed kernel terms correspond to a larger-norm covariance matrix than the same prior over a single kernel term. This implies greater variation in $Y$, preventing the posterior from collapsing around the ADRF (more on this below). These points characterize RIC in the GP response model and suggest that our priors mitigate it.

We note that PRBF performs well in the homogeneous treatment effect experiments relative to A-PRBF, but that this difference is less pronounced or reversed in heterogeneous $\varphi$ settings. This finding agrees with [15]'s argument that independent regularization of $\pi$ and $t$ is necessary for uncovering heterogeneous effects. Specifically, regularizing priors over A-PRBF's summed kernel terms prefer larger-normed GP covariance matrix than the PRBF prior would; roughly speaking, we sum priors over two scale parameters ($\gamma$ and $\sigma$) rather than one ($\gamma$). In our experiments, this assumption is appropriate when we expect greater variation in the response surface (e.g., as a result of more variation in treatment effect), but less so when $\mu$ is linear and/or $\varphi$ is homogeneous. This suggests that, in practical situations where less heterogeneity in treatment effects is expected or the overall response smoother is a relatively smooth process, we should prefer PRBF by default, or look to adjust the scale term priors to reflect that belief.

We observe that non-GP methods with stronger linearity/smoothness priors (i.e., ADS, HI) perform well when treatment effects are homogenous. NPM performs well in the heterogeneous setting, presumably because its product kernel formulation allows it to adapt to a more-rapidly varying response surface and distinguish variation along its $\hat{\pi}$ relative to $t$. BART achieves its best results with $n=500$—in these experiments, we observe a decrease in ADRF bias when the PS is included, but less dramatic than that reported in [15].

**D.3 – IHDP Simulation**

We simulate a response variable using a real continuous treatment and real covariates from the Infant Health and Development Program study [14]. Specifically, we use the number of days in the Child Development Center (CDC) as a continuous treatment. We focus specifically on children with non-zero CDC days. We test two modified configurations of each of the three simulated response surfaces described in [16]:

$$y_i \sim \mathcal{N}(\mu(\mathbf{x}_i) + \eta_i \varphi_i), \text{ where}$$

$$\mu(\mathbf{X}) = \begin{cases} X\beta^A + 4, & A \\ X\beta^B - \omega_B^s, & B \\ Q\beta^C + \omega_C^s, & C \end{cases},$$

$$\eta_i = \begin{cases} 3, & homogeneous \\ 1 + \frac{1}{3}x_{i,1}x_{i,2} + \frac{1}{2}x_{i,1}(1 - x_{i,2}) + 2x_{i,3}, & heterogeneous \end{cases},$$

$$\beta^A = argmax(\mathbf{Z}^A) - 1, \quad \mathbf{z}^A \sim Multinoulli((.5,.2,.15,.1,.05))$$

$$\beta^B = \frac{argmax(\mathbf{Z}^B) - 1}{10}, \quad \mathbf{z}^B_{Cont} \sim Multinoulli((.5,.125,.125,.125,.125))$$

$$\mathbf{z}^B_{Bin} \sim Multinoulli((.6,.1,.1,.1,.1))$$

$$\beta^C = (0, 2.5, 5)_{argmax(\mathbf{z}^C)}, \quad \mathbf{z}^C \sim Multinoulli((.6,.3,.1))$$

where $\varphi_i$ is individual $i$'s CDC days, $x_{i,1}$ is their birthweight, $x_{i,2}$ is their gender label (with 1 indicating female), and $x_{i,3}$ is a composite measure of neo-natal health described in the original study. Note coefficients $\beta$ are as sampled in [16]. We abuse notation to let $\mathbf{Z}^A$ represent the matrix formed by vectors $\mathbf{z}^A$ observed for each covariate and the *argmax* returns a value for each feature. $\mathbf{z}^B_{Cont}$ and $\mathbf{z}^B_{Cont}$ are drawn for each binary and continuous feature, respectively. $(0, 2.5, 5)_{argmax(\mathbf{z}^C)}$ returns the element $argmax(\mathbf{Z}^C)$-th element of the tuple. Finally, $Q$ is the covariate matrix modified to include squared terms for all continuous features and pairwise interactions of all covariates.

Table D.3: IHDP Simulation Results

| Surface | Method | Homogeneous effect | | | | Heterogeneous effect | | | |
|---|---|---|---|---|---|---|---|---|---|
| | | $Cov_{90}$ | $I_{90}$ | Bias | RMSE | $Cov_{90}$ | $I_{90}$ | Bias | RMSE |
| A | A-PRBF | **0.95** | 8.35 | -0.352 | 3.77 | **0.96** | 12.13 | 0.035 | 3.13 |
| | A-RBF | 0.94 | 8.41 | -0.356 | 3.77 | **0.96** | 12.88 | 0.028 | 3.11 |
| | PRBF | **0.95** | 8.35 | -0.336 | 3.74 | **0.96** | 12.72 | 0.048 | 3.10 |
| | RBF | 0.91 | 11.48 | -0.803 | 3.64 | 0.94 | 11.68 | 0.288 | 2.50 |
| | RBF-NT | 0.91 | 11.45 | -0.819 | 3.66 | 0.93 | 11.70 | 0.290 | 2.50 |
| | A-SymG | 0.85 | 19.10 | -0.441 | **3.59** | 0.77 | 14.29 | 0.274 | 3.91 |
| | ADS | 0.47 | 5.68 | -1.610 | 4.74 | 0.50 | 5.13 | -1.519 | 3.82 |
| | BART | - | - | -2.076 | 4.33 | - | - | -1.674 | 3.86 |
| | BART-PS-GBM | - | - | -2.070 | 4.27 | - | - | -1.415 | 3.82 |
| | BART-PS-OLS | - | - | -2.082 | 4.32 | - | - | -1.608 | 3.86 |
| | HI | 0.41 | 4.71 | -1.608 | 4.40 | 0.40 | 4.97 | -1.305 | 3.70 |
| | NW | 0.28 | 3.88 | -2.796 | 5.24 | 0.34 | 4.24 | -1.447 | 3.81 |
| | NPM | **0.66** | 10.01 | -2.038 | 4.70 | **0.68** | 7.67 | -1.776 | 3.94 |
| | IPTW | 0.41 | 5.30 | -1.152 | 3.93 | 0.34 | 3.91 | -0.963 | 3.65 |
| B | A-PRBF | **1.00** | 4.40 | 0.100 | 0.65 | **0.96** | 7.96 | -0.292 | 1.04 |
| | A-RBF | **1.00** | 4.41 | 0.100 | 0.65 | **0.97** | 8.37 | -0.283 | 1.05 |
| | PRBF | **1.00** | 4.71 | 0.038 | 0.60 | **0.97** | 8.47 | -0.252 | 1.13 |
| | RBF | 1.00 | 5.78 | -0.176 | 0.82 | 0.94 | 9.84 | -0.437 | 1.28 |
| | RBF-NT | 1.00 | 5.75 | -0.173 | 0.81 | 0.94 | 9.86 | -0.443 | 1.30 |
| | A-SymG | 1.00 | 11.19 | -0.753 | 1.64 | 0.96 | 8.80 | -0.400 | 1.48 |
| | ADS | 0.40 | 1.10 | -0.001 | 0.58 | 0.53 | 1.51 | 0.150 | 0.98 |
| | BART | - | - | -0.102 | 0.54 | - | - | 0.112 | 0.96 |
| | BART-PS-GBM | - | - | -0.111 | 0.50 | - | - | 0.083 | 0.98 |
| | BART-PS-OLS | - | - | -0.116 | 0.51 | - | - | 0.067 | 0.91 |
| | HI | 0.34 | 1.00 | -0.008 | 0.55 | 0.47 | 1.47 | 0.129 | 0.99 |
| | NW | 0.09 | 1.86 | -1.110 | 2.20 | 0.37 | 1.25 | -0.151 | 1.05 |
| | NPM | **0.68** | 1.46 | 0.037 | 0.59 | **0.72** | 2.33 | -0.135 | 0.89 |
| | IPTW | 0.37 | 0.96 | -0.077 | 0.55 | 0.40 | 1.13 | 0.208 | 1.06 |
| C | A-PRBF | **0.60** | 62.37 | -12.08 | 54.09 | **0.54** | 52.65 | -19.95 | 70.58 |
| | A-RBF | 0.60 | 62.80 | -11.66 | 53.57 | 0.52 | 51.18 | -19.56 | 69.51 |
| | PRBF | 0.57 | 62.64 | -14.56 | 62.30 | 0.53 | 55.86 | -14.13 | 82.71 |
| | RBF | 0.43 | 51.22 | -7.69 | 52.15 | 0.53 | 39.05 | -6.429 | 53.13 |
| | RBF-NT | 0.43 | 51.31 | -7.65 | 53.14 | 0.53 | 39.16 | -6.431 | 53.12 |
| | A-SymG | 0.35 | 49.88 | -12.06 | 55.31 | 0.23 | 48.62 | -11.41 | 54.58 |
| | ADS | 0.48 | 79.50 | -22.56 | 50.34 | 0.30 | 83.11 | -14.32 | 63.24 |
| | BART | - | - | -23.14 | 48.41 | - | - | -8.144 | 58.73 |
| | BART-PS-GBM | - | - | -23.03 | 48.07 | - | - | -8.186 | 59.07 |

| | | | | | | | | |
|---|---|---|---|---|---|---|---|---|
| BART-PS-OLS | - | - | 23.40 | 48.65 | - | - | -8.451 | 58.51 |
| HI | 0.44 | 69.54 | -22.80 | 48.62 | 0.27 | 71.54 | -15.45 | 61.51 |
| NW | 0.31 | 81.45 | -30.19 | 54.27 | 0.16 | 49.11 | -16.39 | 61.00 |
| NPM | **0.69** | **104.8** | -17.53 | 46.43 | 0.38 | 112.1 | -5.880 | 58.79 |
| IPTW | 0.41 | 67.85 | -29.58 | 54.30 | 0.26 | 66.53 | -20.42 | 65.16 |

These experiments reinforce the findings of the study described in D.2, so we summarize:

- PRBF does well (relative to other GP kernels) in simpler processes with homogeneous effects
    - o A-PRBF performs better in more complex processes with heterogeneous effects
- A-PRBF tends to have higher coverage with smaller average interval lengths relative to A-RBF
- RBF and RBF-NT fail to capture all of the variation in ADRF (or in the case of surface B, do so at the expense of much larger mean interval lengths and RMSEs)
- The SymG kernel often greatly under-/over-estimates posterior interval length and tended to have much lower coverage
- ADS and HI do well under homogenous dosage effects, whereas NPM does best in the presence of heterogeneity